\input amstex
\input amsppt.sty
\magnification=\magstep1
\hsize=33.5truecc
\vsize=23truecm
\baselineskip=16truept
\NoBlackBoxes
\TagsOnRight \pageno=1 \nologo
\def\Z{\Bbb Z}
\def\N{\Bbb N}

\def\Q{\Bbb Q}

\def\l{\left}
\def\r{\right}
\def\bg{\bigg}
\def\({\bg(}
\def\[{\bg\lfloor}
\def\){\bg)}
\def\]{\bg\rfloor}
\def\t{\text}
\def\f{\frac}

\def\bi{\binom}
\def\eq{\equiv}

\def\ls{\leqslant}
\def\gs{\geqslant}

\def\M#1#2{\thickfracwithdelims[]\thickness0{#1}{#2}_q}

\def\Proof{\noindent{\it Proof}}

\def\Remark{\medskip\noindent{\it  Remark}}

\def\Ack{\medskip\noindent {\bf Acknowledgments}}
\hbox {Colloq. Math. 154(2018), no.\,2, 241--273.}
\bigskip
\topmatter
\title Two new kinds of numbers and related divisibility results\endtitle
\author Zhi-Wei Sun\endauthor
\leftheadtext{Zhi-Wei Sun}
\affil Department of Mathematics, Nanjing University\\
 Nanjing 210093, People's Republic of China
  \\  zwsun\@nju.edu.cn
  \\ {\tt http://maths.nju.edu.cn/$\sim$zwsun}
\endaffil
\abstract We mainly introduce two new kinds of numbers given by
$$R_n=\sum_{k=0}^n\binom nk\binom{n+k}k\frac1{2k-1}\quad\ (n=0,1,2,\ldots)$$
and
$$S_n=\sum_{k=0}^n\binom nk^2\bi{2k}k(2k+1)\quad\ (n=0,1,2,\ldots).$$
We find that such numbers have many interesting arithmetic properties. For example, if $p\equiv1\pmod 4$ is a prime with $p=x^2+y^2$ (where $x\equiv1\pmod 4$ and $y\equiv0\pmod 2$), then
$$R_{(p-1)/2}\eq p-(-1)^{(p-1)/4}2x\pmod{p^2}.$$ Also,
$$\f1{n^2}\sum_{k=0}^{n-1}S_k\in\Z\ \ \t{and}\ \ \f1n\sum_{k=0}^{n-1}S_k(x)\in\Z[x]\quad\t{for all}\ n=1,2,3,\ldots,$$
where $S_k(x)=\sum_{j=0}^k\bi kj^2\bi{2j}j(2j+1)x^j$. For any positive integers $a$ and $n$, we show that,
somewhat surprisingly,
$$\f1{n^2}\sum_{k=0}^{n-1}(2k+1)\bi{n-1}k^a\bi{-n-1}k^a\in\Z\ \ \t{and}\ \ \f1n\sum_{k=0}^{n-1}\f{\bi{n-1}k^a\bi{-n-1}k^a}{4k^2-1}\in\Z.$$
We also solve a conjecture of V.J.W. Guo and J. Zeng, and pose several conjectures for further research.
\endabstract
\thanks 2010 {\it Mathematics Subject Classification}. \,Primary 11A07, 11B65;
Secondary  05A10, 05A30, 11B75, 11E25.
\newline\indent Supported by the National Natural Science
Foundation (grant 11571162) of China and the NSFC-RFBR Cooperation and Exchange Program (grant 11811530072).
\endthanks

\endtopmatter
\document

\heading{1. Introduction}\endheading

 In combinatorics, the (large) Schr\"oder numbers are given by
 $$S(n)=\sum_{k=0}^n\bi nk\bi{n+k}k\f1{k+1}=\sum_{k=0}^n\bi{n+k}{2k}\bi{2k}k\f1{k+1}\ \ (n\in\N),\tag1.1$$
 where $\N=\{0,1,2,\ldots\}$.
 They are integers since
 $$C_k=\f1{k+1}\bi{2k}k=\bi{2k}k-\bi{2k}{k+1}\in\Z\quad\t{for all}\ k\in\N.$$
 Those $C_n$ with $n\in\N$ are the well-known Catalan numbers. Both Catalan numbers and Schr\"oder numbers have many combinatorial interpretations.
 For example, $S(n)$ is the number of lattice paths from the point $(0,0)$ to $(n,n)$ with only allowed steps $(1,0)$, $(0,1)$ and $(1,1)$ which never rise above the line $y=x$.

 We note that $(2k-1)\mid \bi{2k}k$ for all $k\in\N$. This is obvious for $k=0$. For each $k\in\Z^+=\{1,2,3,\ldots\}$, we have
 $$\f{\bi{2k}k}{2k-1}=\f{2}{2k-1}\bi{2k-1}k=\f 2k\bi{2k-2}{k-1}=2C_{k-1}.$$
 Motivated by this and (1.1), we introduce a new kind of numbers:
 $$R_n:=\sum_{k=0}^n\bi{n}k\bi{n+k}k\f1{2k-1}=\sum_{k=0}^n\bi{n+k}{2k}\bi{2k}k\f1{2k-1}\ \ (n\in\N).\tag1.2$$
Below are the values of $R_0,R_1,\ldots,R_{16}$ respectively:
$$\align&-1,\ 1,\ 7,\ 25,\ 87,\ 329,\ 1359,\ 6001,\ 27759,\ 132689, \ 649815,
\\&3242377,\ 16421831,\ 84196761,\ 436129183,\ 2278835681,\ 11996748255.
\endalign$$
Applying the Zeilberger algorithm (cf. [PWZ, pp.\,101-119]) via {\tt Mathematica 9}, we get the following third-order recurrence for
the new sequence $(R_n)_{n\gs0}$:
$$(n+1)R_n-(7n+15)R_{n+1}+(7n+13)R_{n+2}-(n+3)R_{n+3}=0\quad \t{for}\ n\in\N.\tag1.3$$
In contrast, there is a second-order recurrence for Schr\"oder numbers:
$$nS(n)-3(2n+3)S(n+1)+(n+3)S(n+2)=0\ \ (n=0,1,2,\ldots).$$
So the sequence $(R_n)_{n\gs0}$ looks more sophisticated than Schr\"oder numbers.

For convenience, we also introduce the associated polynomials
$$R_n(x):=\sum_{k=0}^n\bi nk\bi{n+k}k\f{x^k}{2k-1}=\sum_{k=0}^n\bi{n+k}{2k}\bi{2k}k\f{x^k}{2k-1}\in\Z[x].\tag1.4$$
Note that $R_n=R_n(1)$ and $R_n(0)=-1$. Now we list $R_0(x),\ldots,R_5(x)$:
$$\align &R_0(x)=-1,\ R_1(x)=2x-1,\ R_2(x)=2x^2+6x-1,
\\&R_3(x)=4x^3+10x^2+12x-1,\ R_4(x)=10x^4+28x^3+30x^2+20x-1,
\\& R_5(x)=28x^5+90x^4+112x^3+70x^2+30x-1.
\endalign$$
Applying the Zeilberger algorithm via {\tt Mathematica 9}, we get the following third-order recurrence for
the polynomial sequence $(R_n(x))_{n\gs0}$:
$$\aligned&(n+1)R_n(x)-(4nx+10x+3n+5)R_{n+1}(x)+(4nx+6x+3n+7)R_{n+2}(x)
\\&\qquad\qquad\qquad\qquad\qquad\qquad=(n+3)R_{n+3}(x).\endaligned\tag1.5$$

Let $p\eq1\pmod 4$ be a prime. It is well-known that $p$ can be written uniquely as a sum of two squares.
Write $p=x^2+y^2$ with $x\eq1\pmod 4$ and $y\eq0\pmod 2$. In 1828 Gauss (cf. [BEW, (9.0.1)]) proved that
$$\bi{(p-1)/2}{(p-1)/4}\eq2x\pmod p;$$
in 1986 Chowla, Dwork and Evans [CDE] showed further that
$$\bi{(p-1)/2}{(p-1)/4}\eq\f{2^{p-1}+1}2\l(2x-\f p{2x}\r)\pmod{p^2}.$$
The key motivation to introduce the polynomials $R_n(x)\ (x\in\N)$ is our following result.

\proclaim{Theorem 1.1} {\rm (i)} Let $p\eq1\pmod 4$ be a prime, and write $p=x^2+y^2$ with $x\eq1\pmod 4$ and $y\eq0\pmod 2$. Then
$$R_{(p-1)/2}-p\eq\sum_{k=0}^{p-1}\f{\bi{2k}k^2}{(2k-1)(-16)^k}\eq-2\l(\f 2p\r)x\pmod{p^2},\tag1.6$$
where $(\f{\cdot}p)$ denotes the Legendre symbol.
Also,
$$R_{(p-1)/2}(-2)+2p\l(\f 2p\r)\eq\sum_{k=0}^{p-1}\f{\bi{2k}k^2}{(2k-1)8^k}\eq\l(\f2p\r)\f p{2x}\pmod{p^2},\tag1.7$$
$$R_{(p-1)/2}\l(-\f12\r)+\f p2\l(\f 2p\r)\eq\sum_{k=0}^{p-1}\f{\bi{2k}k^2}{(2k-1)32^k}\eq\f p{4x}-x\pmod{p^2}.\tag1.8$$

{\rm (ii)} Let $p\eq3\pmod 4$ be a prime. Then
$$R_{(p-1)/2}\eq \sum_{k=0}^{p-1}\f{\bi{2k}k^2}{(2k-1)(-16)^k}\eq-\f12\l(\f 2p\r)\bi{(p+1)/2}{(p+1)/4}\pmod p\tag1.9$$
and
$$R_{(p-1)/2}(-2)\eq \sum_{k=0}^{p-1}\f{\bi{2k}k^2}{(2k-1)8^k}\eq-\f12\l(\f 2p\r)\bi{(p+1)/2}{(p+1)/4}\pmod p\tag1.10$$
We also have
$$R_{(p-1)/2}\l(-\f12\r)+\f p2\l(\f2p\r)\eq\sum_{k=0}^{p-1}\f{\bi{2k}k^2}{(2k-1)32^k}\eq-\f{p+1}{2^p+2}\bi{(p+1)/2}{(p+1)/4}\pmod {p^2}.\tag1.11$$
\endproclaim

Our following theorem is motivated by (1.7).

\proclaim{Theorem 1.2} Let $p=2n+1$ be any odd prime. Then
$$\sum_{k=0}^{p-1}\f{\bi{2k}k\bi{2k}{k+d}}{(2k-1)8^k}\eq0\pmod p\tag1.12$$
for all $d\in\{0,\ldots,n\}$ with $d\eq n\pmod 2$.
\endproclaim
\Remark\ 1.1. In contrast with (1.12), by induction we have
$$\sum_{k=0}^n\f{\bi{2k}k\bi{2k}{k+d}}{(2k-1)16^k}=\f{2n+1}{(4d^2-1)16^n}\bi{2n}n\bi{2n}{n+d}\quad\t{for all}\ d,n\in\N.$$

Below is our third theorem.

\proclaim{Theorem 1.3} {\rm (i)} For any odd prime $p$, we have
$$\sum_{k=0}^{p-1}R_k\eq-p-\l(\f{-1}p\r)\pmod{p^2}.\tag1.13$$

{\rm (ii)} For any positive integer $n$, we have
$$R_n(-1)=-(2n+1)\tag1.14$$
and consequently
$$\sum_{k=0}^n\f{\bi nk\bi{-n}k}{2k-1}=-2n.\tag1.15$$
\endproclaim
\Remark\ 1.2. Although there are many known combinatorial identities (cf. [G]), (1.15) seems new and concise.
\medskip

Now we introduce another kind of new numbers:
$$S_n:=\sum_{k=0}^n\bi nk^2\bi{2k}k(2k+1)\quad\ (n=0,1,2,\ldots).\tag1.16$$
We also define the associated polynomials
$$S_n(x):=\sum_{k=0}^n\bi nk^2\bi{2k}k(2k+1)x^k\quad (n=0,1,2,\ldots).\tag1.17$$

Here are the values of $S_0,S_1,\ldots,S_{12}$ respectively:
$$\align&1,\ 7,\ 55,\ 465,\ 4047,\ 35673,\ 316521,\ 2819295,\ 25173855,
\\&225157881,\ 2016242265,\ 18070920255,\ 162071863425.
\endalign$$
Now we list the polynomials $S_0(x),\ldots,S_5(x)$:
$$\align &S_0(x)=1,\ S_1(x)=6x+1,\ S_2(x)=30x^2+24x+1,
\\&S_3(x)=140x^3+270x^2+54x+1,
\\&S_4(x)=630x^4+2240x^3+1080x^2+96x+1,
\\& S_5(x)=2772x^5+15750x^4+14000x^3+3000x^2+150x+1.
\endalign$$
Applying the Zeilberger algorithm via {\tt Mathematica 9}, we get the following recurrence for $(S_n)_{n\gs0}$:
$$9(n+1)^2S_n-(19n^2+74n+87)S_{n+1}+(n+3)(11n+29)S_{n+2}=(n+3)^2S_{n+3},\tag1.18$$
which looks more complicated than the recurrence relation (1.3) for $(R_n)_{n\gs0}$.
Also, the Zeilberger algorithm could yield a very complicated third-order recurrence for
the polynomial sequence $(S_n(x))_{n\gs0}$. Despite these complicated recurrences, we are able to establish the following result which looks interesting.

\proclaim{Theorem 1.4} {\rm (i)} For any positive integer $n$, we have
$$\f1{n^2}\sum_{k=0}^{n-1}S_k=\sum_{k=0}^{n-1}\bi{n-1}k^2C_k\in\Z\tag1.19$$
and
$$\f1n\sum_{k=0}^{n-1}S_k(x)\in\Z[x].\tag1.20$$

{\rm (ii)} For any prime $p>3$, we have
$$\sum_{k=1}^{p-1}\f{S_k}k\eq p\sum_{k=1}^{p-1}\f{S_k}{k^2}\eq-\f p2\l(\f p3\r)B_{p-2}\l(\f13\r)\pmod{p^2},$$
where $B_n(x)$ denotes the Bernoulli polynomial of degree $n$.
\endproclaim

In 2012 Guo and Zeng [GZ, Corollary 5.6] employed $q$-binomial coefficients to prove that for any $a,b\in\N$ and positive integer $n$ we have
$$\sum_{k=0}^{n-1}(-1)^{(a+b)k}\bi{n-1}k^a\bi{-n-1}k^b\eq0\pmod n.$$
(Note that $\bi{-n-1}k=(-1)^k\bi{n+k}k$.)
This, together with (1.15) and Theorem 1.4, led us to obtain the following result via a new method.

\proclaim{Theorem 1.5} {\rm (i)} Let $a_1,\ldots,a_m$ and $n>0$ be integers. Then
$$\align\sum_{k=0}^{n-1}(\pm1)^k(2k+1)\prod_{i=1}^m\bi{a_in-1}k\eq&0\pmod n,\tag1.21
\\\sum_{k=0}^{n-1}(\pm1)^k(4k^3-1)\prod_{i=1}^m\bi{a_in-1}k\eq&0\pmod n.\tag1.22
\endalign$$
Also,
$$\gcd(a_1+\cdots+a_m-1,2)\sum_{k=0}^{n-1}(-1)^{km}(2k+1)\prod_{i=1}^m\bi{a_in-1}k\eq0\pmod {n^2},\tag1.23$$
and
$$6\sum_{k=0}^{n-1}(-1)^{km}(3k^2+3k+1)\prod_{i=1}^m\bi{a_in-1}k\eq0\pmod {n^2}.\tag1.24$$
Moreover,
$$\sum_{k=0}^{n-1}(-1)^k(4k^3-1)\prod_{i=1}^m\bi{a_in-1}k\bi{-a_in-1}k\eq0\pmod{n^2},\tag1.25$$
and
$$\aligned\gcd(a_1+\cdots+a_m-1,2)&\sum_{k=0}^{n-1}(3k^2+3k+1)\prod_{i=1}^m\bi{a_in-1}k\bi{-a_in-1}k
\\&\qquad\qquad\eq0\pmod{n^3}.
\endaligned\tag1.26$$

{\rm (ii)} For any positive integers $a,b,n$, we have
$$\gather \f1n\sum_{k=0}^{n-1}\f{\bi{n-1}k^a\bi{-n-1}k^a}{4k^2-1}\in\Z,\ \ \f1n\sum_{k=0}^{n-1}\f{\bi{n-1}k^a\bi{-n-1}k^a}{\bi{k+2}2}\in\Z,\tag1.27
\\\f1n\sum_{k=0}^{n-1}(-1)^k\l(1+\f{2k}{4k^2-1}\r)\bi{n-1}k^a\bi{-n-1}k^a\in\Z,\tag1.28
\\\f1n\sum_{k=0}^{n-1}(-1)^k\(4-\f{2k+3}{\bi{k+2}2}\)\bi{n-1}k^a\bi{-n-1}k^a\in\Z,\tag1.29
\endgather$$
and
$$\gather\sum_{k=0}^{n-1}\f{(-1)^{(a+b)k}}{4k^2-1}\bi{n-1}k^a\bi{-n-1}k^b\in\Z,\tag1.30
\\\sum_{k=0}^{n-1}\f{(-1)^{(a+b-1)k}k}{4k^2-1}\bi{n-1}k^a\bi{-n-1}k^b\in\Z,\tag1.31
\\\sum_{k=0}^{n-1}\f{(-1)^{(a+b)k}}{\bi{k+2}2}\bi{n-1}k^a\bi{-n-1}k^b\in\Z,\tag1.32
\\\sum_{k=0}^{n-1}\f{(-1)^{(a+b-1)k}(2k+3)}{\bi{k+2}2}\bi{n-1}k^a\bi{-n-1}k^b\in\Z,\tag1.33
\\\sum_{k=0}^{n-1}\f{(-1)^{(a+b)k}(3k+1)}{(2k+1)\bi{2k}k}\bi{n-1}k^a\bi{-n-1}k^b\in\Z,\tag1.34
\\\sum_{k=0}^{n-1}\f{(-1)^{(a+b-1)k}(5k+3)}{(2k+1)\bi{2k}k}\bi{n-1}k^a\bi{-n-1}k^b\in\Z.\tag1.35
\endgather$$
\endproclaim
\Remark\ 1.3. For any positive integer $n$, using (1.15) we can deduce that
$$\sum_{k=1}^{n-1}\f{\bi{n-1}k\bi{-n-1}k}{4k^2-1}=\f12\sum_{k=0}^n\f{\bi nk\bi{-n}k}{2k-1}=-n.$$
An extension of (1.21) given in (4.5) confirms a conjecture of Guo and Zeng [GZ].
By (1.23), for any positive integers $a,b,n$ we have the congruence
$$\gcd(a+b-1,2)\sum_{k=0}^{n-1}(-1)^{(a+b)k}(2k+1)\bi{n-1}k^a\bi{-n-1}k^b\eq0\pmod{n^2}.\tag1.36$$

\proclaim{Corollary 1.1} For $n\in\N$ define
$$\align t_n=&\sum_{k=0}^n\bi nk^2\bi{n+k}k^2\f1{2k-1},
\\T_n=&\sum_{k=0}^n\bi nk^2\bi{n+k}k^2(2k+1),
\\T_n^+=&\sum_{k=0}^n\bi nk^2\bi{n+k}k^2(2k+1)^2,
\\T_n^-=&\sum_{k=0}^n\bi nk^2\bi{n+k}k^2(-1)^k(2k+1)^2.
\endalign$$
Then, for any positive integer $n$, we have
$$\f1{n^3}\sum_{k=0}^{n-1}(2k+1)t_k\in\Z,\ \f1{n^3}\sum_{k=0}^{n-1}(2k+1)T_k\in\Z,\tag1.37$$
and
$$\f1{n^4}\sum_{k=0}^{n-1}(2k+1)T_k^+\in\Z,\ \f1{n^3}\sum_{k=0}^{n-1}(2k+1)T_k^{-}\in\Z.\tag1.38$$
\endproclaim

We will prove Theorems 1.1-1.3 in the next section. We are going to show Theorem 1.4 and a $q$-congruence related to (1.21) in Section 3.
Section 4 is devoted to our proofs of Theorem 1.5 and Corollary 1.1 and some extensions.
In Section 5 we pose several related conjectures for further research.

\heading{2. Proofs of Theorems 1.1-1.3}\endheading

\proclaim{Lemma 2.1} Let $p=2n+1$ be an odd prime. Then
$$\aligned R_n(x)\eq&\sum_{k=0}^{n}\f{\bi{2k}k^2}{2k-1}\l(-\f x{16}\r)^k
\\\eq&\sum_{k=0}^{p-1}\f{\bi{2k}k^2}{2k-1}\l(-\f x{16}\r)^k-p(-x)^{n+1}\pmod{p^2}.
\endaligned\tag2.1$$
\endproclaim
\Proof. As pointed out in [S11, Lemma 2.2], for each $k=0,\ldots,n$ we have
$$\bi{(p-1)/2+k}{2k}=\f{\prod_{0<j\ls k}(p^2-(2j-1)^2)}{(2k)!4^k}\eq\f{\bi{2k}k}{(-16)^k}\pmod{p^2}.$$
Recall that $(2k-1)\mid\bi{2k}k$ for all $k\in\N$. Therefore,
$$R_n(x)=\sum_{k=0}^n\bi{n+k}{2k}\bi{2k}k\f{x^k}{2k-1}\eq\sum_{k=0}^n\f{\bi{2k}k^2}{2k-1}\l(-\f{x}{16}\r)^k\pmod{p^2}.$$
Clearly, $p\mid\bi{2k}k$ for all $k=n+1,\ldots,p-1$. Also,
$$\align\f{\bi{p+1}{(p+1)/2}^2}{2\times(p+1)/2-1}\l(-\f x{16}\r)^{(p+1)/2}
=&\f{4p\bi{p-1}{(p-3)/2}^2}{((p-1)/2)^2}\times\f{(-x)^{(p+1)/2}}{4^{p+1}}
\\\eq& p(-x)^{(p+1)/2}\pmod{p^2}.
\endalign$$
So the second congruence in (2.1) also holds. \qed

\proclaim{Lemma 2.2} For any nonnegative integer $n$, we have
$$\sum_{k=0}^n\l((16-x)k^2-4\r)\f{\bi{2k}k^2}{2k-1}x^{n-k}=\f{4(n+1)^2}{2n+1}\bi{2n+1}n^2.\tag2.2$$
\endproclaim
\Proof. Let $P(x)$ denote the left-hand side of (2.2). Then
$$\align P(x)=&4\sum_{k=0}^n(4k^2-1)\f{\bi{2k}k^2}{2k-1}x^{n-k}-\sum_{k=0}^n\f{k^2\bi{2k}k^2}{2k-1}x^{n+1-k}
\\=&4\sum_{k=0}^n(2k+1)\bi{2k}k^2x^{n-k}-4\sum_{k=1}^n(2k-1)\bi{2(k-1)}{k-1}^2x^{n-(k-1)}
\\=&4(2n+1)\bi{2n}n^2=\f{4(n+1)^2}{2n+1}\bi{2n+1}{n+1}^2.
\endalign$$
This concludes the proof. \qed

\medskip
\noindent{\it Proof of Theorem} 1.1. Applying Lemma 2.1 with $x=1,-2,-1/2$ we get the first congruence in each of (1.6)-(1.11).

Let $p$ be an odd prime. For any $p$-adic integer $m\not\eq0\pmod p$, by Lemma 2.2 we have
$$(16-m)\sum_{k=1}^{p-1}\f{k^2\bi{2k}k^2}{(2k-1)m^k}-4\sum_{k=0}^{p-1}\f{\bi{2k}k^2}{(2k-1)m^k}\eq0\pmod {p^2}$$
and hence
$$\align\sum_{k=0}^{p-1}\f{\bi{2k}k^2}{(2k-1)m^k}\eq&(16-m)\sum_{k=1}^{p-1}(2k-1)\f{\bi{2(k-1)}{k-1}^2}{m^k}
\\=&\l(\f{16}m-1\r)\(\sum_{j=0}^{p-1}(2j+1)\f{\bi{2j}j^2}{m^j}-(2p-1)\f{\bi{2p-2}{p-1}^2}{m^{p-1}}\)
\\\eq&\l(\f{16}m-1\r)\sum_{k=0}^{(p-1)/2}(2k+1)\f{\bi{2k}k^2}{m^k}\pmod{p^2}.
\endalign$$ (Note that $\bi{2p-2}{p-1}=(2p-2)!/((p-1)!)^2\eq0\pmod p$.)
Taking $m=-16,8,32$ we obtain
$$\align\sum_{k=0}^{p-1}\f{\bi{2k}k^2}{(2k-1)(-16)^k}\eq&-2\sum_{k=0}^{(p-1)/2}(2k+1)\f{\bi{2k}k^2}{(-16)^k}\pmod{p^2},\tag2.3
\\\sum_{k=0}^{p-1}\f{\bi{2k}k^2}{(2k-1)8^k}\eq&\sum_{k=0}^{(p-1)/2}(2k+1)\f{\bi{2k}k^2}{8^k}\pmod{p^2},\tag2.4
\\\sum_{k=0}^{p-1}\f{\bi{2k}k^2}{(2k-1)32^k}\eq&-\f12\sum_{k=0}^{(p-1)/2}(2k+1)\f{\bi{2k}k^2}{32^k}\pmod{p^2}.\tag2.5
\endalign$$

(i) Recall the condition $p=x^2+y^2$ with $x\eq1\pmod 4$ and $y\eq0\pmod 2$. By [Su12a, Theorem 1.2],
$$\l(\f2p\r)x\eq\sum_{k=0}^{(p-1)/2}\f{2k+1}{(-16)^k}\bi{2k}k^2\eq\sum_{k=0}^{(p-1)/2}\f{k+1}{8^k}\bi{2k}k^2\pmod{p^2}.$$
The author [Su11, Conjecture 5.5] conjectured that
$$\l(\f 2p\r)\sum_{k=0}^{(p-1)/2}\f{\bi{2k}k^2}{8^k}\eq\sum_{k=0}^{(p-1)/2}\f{\bi{2k}k^2}{32^k}\eq2x-\f p{2x}\pmod{p^2}$$
which was later confirmed by the author's brother Z.-H. Sun [S11], who also showed that
$$\sum_{k=0}^{(p-1)/2}\f{k\bi{2k}k^2}{32^k}\eq0\pmod{p^2}.$$
Combining these with (2.3)-(2.5), we immediately get the second congruences in (1.6)-(1.8).

(ii) Now we consider the case $p\eq3\pmod 4$. By [Su13a, Theorem 1.3],
$$\sum_{k=0}^{(p-1)/2}\f{2k+1}{(-16)^k}\bi{2k}k^2\eq\sum_{k=0}^{(p-1)/2}\f{2k}{(-16)^k}\bi{2k}k^2\eq\f14\l(\f2p\r)\bi{(p+1)/2}{(p+1)/4}\pmod{p}$$
and
$$\sum_{k=0}^{(p-1)/2}\f{2k+1}{8^k}\bi{2k}k^2\eq\sum_{k=0}^{(p-1)/2}\f{2k}{8^k}\bi{2k}k^2\eq-\f12\l(\f2p\r)\bi{(p+1)/2}{(p+1)/4}\pmod{p}.$$
Combining this with (2.3) and (2.4), we obtain the second congruences in (1.9) and (1.10).

Z.-H. Sun [S11, Theorem 2.2] confirmed the author's conjectural congruence
$$\sum_{k=0}^{(p-1)/2}\f{\bi{2k}k^2}{32^k}\eq0\pmod{p^2}.$$
He also showed [S11, Theorem 2.3] that
$$\sum_{k=0}^{(p-1)/2}\f{k\bi{2k}k^2}{32^k}\eq\l(\f 2p\r)\f{p+1}{4\times 2^{(p-1)/2}}\bi{(p+1)/2}{(p+1)/4}\pmod{p^2}.$$
Observe that
$$\align2^{p-1}+1=&2+\l(\l(\f2 p\r)2^{(p-1)/2}+1\r)\l(\l(\f2 p\r)2^{(p-1)/2}-1\r)
\\\eq&2+2\l(\l(\f2 p\r)2^{(p-1)/2}-1\r)=2\l(\f2 p\r)2^{(p-1)/2}\pmod{p^2}.
\endalign$$
Therefore,
$$\sum_{k=0}^{(p-1)/2}\f{2k+1}{32^k}\bi{2k}k^2\eq\f{p+1}{2^{p-1}+1}\bi{(p+1)/2}{(p+1)/4}\pmod{p^2}.$$
Combining this with (2.5) we obtain the second congruence in (1.11).

The proof of Theorem 1.1 is now complete. \qed

\medskip\noindent
{\it Proof of Theorem} 1.2. Clearly $\bi{2k}k/(2k-1)\eq0\pmod p$ if $n+1<k<p$. Thus
$$\align\sum_{k=0}^{p-1}\f{\bi{2k}k\bi{2k}{k+n}}{(2k-1)8^k}
\eq&\sum_{k=n}^{n+1}\f{\bi{2k}k\bi{2k}{k+n}}{(2k-1)8^k}
\\=&\f{\bi{p-1}n}{(2n-1)8^n}+\f{\bi{p+1}{n+1}\bi{p+1}p}{p8^{n+1}}
\\\eq&\f12(-1)^{n+1}\l(\f 8p\r)+\f{2\f pn\bi{p-1}{n-1}(p+1)}{p8^{n+1}}\eq0\pmod p.
\endalign$$
So (1.12) holds for $d=n$.

 Define
 $$u_m(d)=\sum_{k=0}^m\f{\bi{2k}k\bi{2k}{k+d}}{(2k-1)8^k}\quad\t{for}\ d,m\in\N.$$
Applying the Zeilberger algorithm via {\tt Mathematica 9}, we get the recurrence
$$(2d-1)u_m(d)+(2d+5)u_m(d+2)=(d+1)\f{\bi{2m}m\bi{2m+2}{m+d+2}}{(m+1)8^m}.$$
If $0\ls d\ls n-2$, then
$$\f{\bi{2(p-1)}{p-1}\bi{2p}{p+d+1}}{8^{p-1}p}=\f {\f p{2p-1}\bi{2p-1}p\f{2p}{p+d+1}\bi{2p-1}{p+d}}{8^{p-1}p}\eq0\pmod p$$
and hence
$$(2d-1)u_{p-1}(d)\eq-(2d+5)u_{p-1}(d+2)\pmod p,$$
therefore
$$u_{p-1}(d+2)\eq0\pmod p\ \Longrightarrow\ u_{p-1}(d)\eq0\pmod p.$$

In view of the above, we have proved the desired result by induction. \qed

\proclaim{Lemma 2.3} For any integers $k>0$ and $n\gs0$, we have the identity
$$(-1)^k\f{\bi nk\bi{-n}k}{\bi{2k-1}k}=\f{2n}{n+k}\bi{n+k}{2k}=\bi{n+k}{2k}+\bi{n+k-1}{2k}.\tag2.6$$
\endproclaim
\Proof. Observe that
$$\align(-1)^k\bi nk\bi{-n}k=&\bi nk\bi{n+k-1}k=\bi nk\bi{n+k}k\f n{n+k}
\\=&\bi{n+k}{2k}\bi{2k}k\f n{n+k}=\f{2n}{n+k}\bi{n+k}{2k}\bi{2k-1}k
\endalign$$
and
$$\f{2n}{n+k}\bi{n+k}{2k}=\l(1+\f{n-k}{n+k}\r)\bi{n+k}{2k}=\bi{n+k}{2k}+\bi{n+k-1}{2k}.$$
So (2.6) follows. \qed

\medskip\noindent{\it Proof of Theorem} 1.3. (i) It is known that
$$\sum_{n=0}^m\bi{n+l}l=\bi{l+m+1}{l+1}\quad\t{for all}\ \ l,m\in\N$$
(cf. [G, (1.49)]). Thus
$$\align\sum_{n=0}^{p-1}R_n=&\sum_{n=0}^{p-1}\sum_{k=0}^n\bi{n+k}{2k}\f{\bi{2k}k}{2k-1}=\sum_{k=0}^{p-1}\f{\bi{2k}k}{2k-1}\sum_{n=k}^{p-1}\bi{n+k}{2k}
\\=&\sum_{k=0}^{p-1}\f{\bi{2k}k}{2k-1}\bi{p+k}{2k+1}=\sum_{k=0}^{p-1}\f p{(2k+1)(2k-1)}\prod_{0<j\ls k}\f{p^2-j^2}{j^2}
\\\eq&p\sum_{k=0}^{p-1}\f{(-1)^k}{4k^2-1}=-p+p\sum_{k=1}^{(p-1)/2}\l(\f{(-1)^k}{4k^2-1}+\f{(-1)^{p-k}}{4(p-k)^2-1}\r)
\\\eq&-p+p\l(\f{(-1)^{(p-1)/2}}{4((p-1)/2)^2-1}+\f{(-1)^{(p+1)/2}}{4((p+1)/2)^2-1}\r)
\\\eq&-p+\l(\f{-1}p\r)\l(\f1{p-2}-\f1{p+2}\r)\eq-p-\l(\f{-1}p\r)\pmod{p^2}.
\endalign$$

(ii) For any positive integer $n$, clearly
$$\align &R_{n}(-1)-R_{n-1}(-1)
\\=&\sum_{k=0}^{n}\l(\bi{n+k}{2k}-\bi{n-1+k}{2k}\r)\bi{2k}k\f{(-1)^k}{2k-1}
\\=&\sum_{k=1}^{n}\bi{n-1+k}{2k-1}(-1)^k2C_{k-1}
=-2\sum_{j=0}^{n-1}\bi{n+j}{2j+1}(-1)^jC_j
\endalign$$
and hence $R_{n}(-1)-R_{n-1}(-1)=-2$ with the help of [Su12b, (2.6)]. Thus, by induction, (1.14) holds for all $n\in\N$.

In view of (2.6) and (1.14), for each positive integer $n$ we have
$$\align2\sum_{k=0}^n\f{\bi nk\bi{-n}k}{2k-1}=&\sum_{k=0}^n\l(\bi{n+k}{2k}+\bi{n+k-1}{2k}\r)\bi{2k}k\f{(-1)^k}{2k-1}
\\=&\sum_{k=0}^n\l(\bi nk\bi{n+k}k+\bi{n-1}k\bi{n-1+k}k\r)\f{(-1)^k}{2k-1}
\\=&R_n(-1)+R_{n-1}(-1)=-(2n+1)-(2n-1)=-4n
\endalign$$
and hence (1.15) holds.

The proof of Theorem 1.3 is now complete. \qed

\heading{3. Proof of Theorem 1.4 and a $q$-congruence related to (1.21)}\endheading

\medskip
\noindent{\it Proof of} (1.19). Define
$$h_n:=\sum_{k=0}^n\bi nk^2C_k\quad\t{for}\ n=0,1,2,\ldots.$$
We want to show that $\sum_{k=0}^{n-1}S_k=n^2h_{n-1}$ for any positive integer $n$.
This is trivial for $n=1$. So, it suffices to show that
$$S_n=(n+1)^2h_n-n^2h_{n-1}=\sum_{k=0}^n((n+1)^2-(n-k)^2)\bi nk^2C_k$$
for all $n=1,2,3,\ldots$. Define $v_n=\sum_{k=0}^n((n+1)^2-(n-k)^2)\bi nk^2C_k$ for $n\in\N$.
It is easy to check that $v_n=S_n$ for $n=0,1,2$. Via the Zeilberger algorithm we find the recurrence
$$9(n+1)^2v_n-(19n^2+74n+87)v_{n+1}+(n+3)(11n+29)v_{n+2}=(n+3)^2v_{n+3}.$$
This, together with (1.18), implies that $v_n=S_n$ for all $n\in\N$. \qed

For each integer $n$ we set
$$[n]_q=\f{1-q^n}{1-q},$$
which is the usual $q$-analogue of $n$. For any $n\in\Z$, we define
$$\M n0=1\quad\t{and}\quad \M nk=\f{\prod_{j=0}^{k-1}[n-j]_q}{\prod_{j=1}^k[j]_q}\ \ \t{for}\ k=1,2,3,\ldots.$$
Obviously $\lim_{q\to1}\M nk=\bi nk$ for all $k\in\N$ and $n\in\Z$. It is easy to see that
$$\M nk=q^k\M{n-1}k+\M{n-1}{k-1}\quad \ \t{for all}\ k,n=1,2,3,\ldots.$$
By this recursion, $\M nk\in\Z[q]$ for all $k,n\in\N$. For any integers $a,\,b$ and $n>0$, clearly
$$a\eq b\pmod n\ \Longrightarrow\ [a]_q\eq[b]_q\pmod{[n]_q}.$$

Let $n$ be a positive integer. The cyclotomic polynomial
$$\Phi_n(q):=\prod^n\Sb a=1\\(a,n)=1\endSb\l(q-e^{2\pi ia/n}\r)\in\Z[q]$$
is irreducible in the ring $\Z[q]$. It is well-known that
$$q^n-1=\prod_{d\mid n}\Phi_d(q).$$
Note that $\Phi_1(q)=q-1$.

\proclaim{Lemma 3.1 {\rm ($q$-Lucas Theorem (cf. [O]))}} Let $a,b,d,s,t\in\N$ with $s<d$ and $t<d$. Then
$$\M{ad+s}{bd+t}\eq\bi ab\M st\pmod{\Phi_d(q)}.\tag3.1$$
\endproclaim

\proclaim{Lemma 3.2} Let $n$ be a positive integer and let $k\in\N$ with $k<(n-1)/2$. Then
$$\sum_{h=0}^{n-1}q^h\M hk^2\eq0\pmod{\Phi_n(q)}.\tag3.2$$
\endproclaim
\Proof. Note that
$$\sum_{h=0}^{n-1}q^h\M hk^2=\sum_{m=0}^{n-1-k}q^{k+m}\M{k+m}m^2$$
and
$$\align\M{k+m}m=&\prod_{j=1}^m\f{1-q^{k+j}}{1-q^j}=\prod_{j=1}^m\l(q^{k+j}\f{q^{-k-j}-1}{1-q^j}\r)
\\=&(-1)^mq^{km+m(m+1)/2}\prod_{j=1}^m\f{1-q^{-k-j}}{1-q^j}=(-1)^mq^{km+m(m+1)/2}\M{-k-1}m.
\endalign$$
Thus
$$\align\sum_{h=0}^{n-1}q^h\M hk^2=&\sum_{m=0}^{n-1-k}q^{k+m}q^{2km+m(m+1)}\M{-k-1}m^2
\\\eq&q^{-k^2-k-1}\sum_{m=0}^{n-1-k}q^{(k+m+1)^2}\M{n-k-1}m^2
\\\eq&q^{-k(k+1)-1}\sum_{m=0}^{n-1-k}q^{(n-k-m-1)^2}\M{n-k-1}m\M{n-k-1}{n-k-1-m}
\\=&q^{-k(k+1)-1}\M{2(n-k-1)}{n-k-1}\pmod{\Phi_n(q)}
\endalign$$
with the help of the $q$-Chu-Vandermonde identity (cf. [AAR, p.\,542]).
As $2(n-1-k)\gs n>n-k-1$, $\M{2(n-k-1)}{n-k-1}$ is divisible by $\Phi_n(q)$.
Therefore (3.2) holds. \qed

\proclaim{Theorem 3.1} For any integers $n>k\gs0$, we have
$$[2k+1]_q\M{2k}k\sum_{h=0}^{n-1}q^h\M{h}k^2\eq0\pmod{[n]_q}\tag3.3$$
and hence
$$(2k+1)\bi{2k}k\sum_{h=0}^{n-1}\bi hk^2\eq0\pmod n.\tag3.4$$
\endproclaim
\Proof. Clearly (3.3) with $q\to1$ yields (3.4), and (3.3) holds trivially in the case $n=1$ and $k=0$.
Below we only need to prove (3.3) for $n>1$.

As the polynomials $\Phi_2(q),\Phi_3(q),\ldots$ are pairwise coprime and
$$[n]_q=\prod\Sb d\mid n\\d>1\endSb\Phi_d(q),\tag3.5$$
it suffices to show
$$[2k+1]_q\M{2k}k\sum_{h=0}^{n-1}q^h\M{h}k^2\eq0\pmod{\Phi_d(q)}\tag3.6$$
for any divisor $d>1$ of $n$. Set $m=n/d$ and write $k=bd+t$ with $b,t\in\N$ and $t<d$. If $t<(d-1)/2$, then by applying Lemmas 3.1 and 3.2 we obtain
$$\align \sum_{h=0}^{n-1}q^h\M hk^2
=&\sum_{a=0}^{m-1}\sum_{s=0}^{d-1}q^{ad+s}\M{ad+s}{bd+t}^2
\\\eq&\sum_{a=0}^{m-1}\sum_{s=0}^{d-1}q^s\bi ab^2\M{s}{t}^2=\sum_{a=0}^{m-1}\bi ab^2\,\sum_{s=0}^{d-1}q^s\M st^2
\eq0\pmod{\Phi_d(q)}.\endalign$$
If $t=(d-1)/2$, then
$$[2k+1]_q=[2bd+2t+1]_q=[(2b+1)d]_q\eq0\pmod{[d]_q}.$$
When $d/2\ls t<d$, by Lemma 3.1 we have
$$\M{2k}k=\M{(2b+1)d+2t-d}{bd+t}\eq\bi{2b+1}b\M{2t-d}t=0\pmod{\Phi_d(q)}.$$
So (3.6) holds, and this completes the proof. \qed

\medskip
\noindent{\it Proof of} (1.20). In light of (3.4),
$$\align\f 1n\sum_{h=0}^{n-1}S_h(x)=&\f1n\sum_{h=0}^{n-1}\sum_{k=0}^h\bi hk^2\bi{2k}k(2k+1)x^k
\\=&\sum_{k=0}^{n-1}\f{x^k}n(2k+1)\bi{2k}k\sum_{h=0}^{n-1}\bi hk^2\in\Z[x].
\endalign$$
This concludes the proof. \qed

\medskip
\noindent{\it Proof of Theorem} 1.4(ii). Let $p>3$ be a prime. By a well-known result of Wolstenholme [W],
$$\sum_{k=1}^{p-1}\f1k\eq0\pmod{p^2}\quad\t{and}\quad\sum_{k=1}^{p-1}\f1{k^2}\eq0\pmod p.$$
Clearly,
$$\align\sum_{n=1}^{p-1}\f{S_n}{n^2}=&\sum_{n=1}^{p-1}\f1{n^2}\sum_{k=0}^n\bi nk^2\bi{2k}k(2k+1)
\\\eq&\sum_{k=1}^{p-1}\bi{2k}k\f{2k+1}{k^2}\sum_{n=k}^{p-1}\bi{n-1}{k-1}^2
\\=&\sum_{k=1}^{p-1}\f{2k+1}{k^3}2(2k-1)\bi{2(k-1)}{k-1}\sum_{h=0}^{p-1}\bi h{k-1}^2
\\&-\sum_{k=1}^{p-1}\f{2k+1}{k^2}\bi{2k}k\bi{p-1}{k-1}^2
\\\eq&-\sum_{k=1}^{p-1}\f{2k+1}{k^2}\bi{2k}k\pmod{p}
\endalign$$
with the help of Theorem 3.1.
Note that
$$\sum_{k=1}^{p-1}\f{\bi{2k}k}{k}\eq0\pmod{p^2}\quad \t{and}\quad \sum_{k=1}^{p-1}\f{\bi{2k}k}{k^2}\eq\f12\l(\f p3\r)B_{p-2}\l(\f13\r)\pmod p$$
by [ST] and [MT] respectively. Therefore,
$$\sum_{n=1}^{p-1}\f{S_n}{n^2}\eq-2\sum_{k=1}^{p-1}\f{\bi{2k}k}k-\sum_{k=1}^{p-1}\f{\bi{2k}k}{k^2}\eq-\f12\l(\f p3\r)B_{p-2}\l(\f13\r)\pmod p.$$
Observe that
$$\align\sum_{n=1}^{p-1}\f{S_n}n=&\sum_{n=1}^{p-1}\f1n\sum_{k=0}^n\bi nk^2\bi{2k}k(2k+1)
\\\eq&\sum_{k=1}^{p-1}\bi{2k}k\f{2k+1}k\sum_{n=k}^{p-1}\bi{n-1}{k-1}\bi nk
\\\eq&\sum_{k=1}^{p-1}\bi{2k}k\f{2k+1}k\(\sum_{n=k}^{p-1+k}\bi{n-1}{k-1}\bi nk-\bi{p-1}{k-1}\bi pk\)
\\=&\sum_{k=1}^{p-1}\f{2k+1}k\bi{2k}k\(\sum_{j=0}^{p-1}\bi{k+j-1}j\bi{k+j}j-\f pk\bi{p-1}{k-1}^2\)
\\\eq&\sum_{k=1}^{p-1}\f{2k+1}k\bi{2k}k\sum_{j=0}^{p-1}\bi{-k}j\bi{-k-1}j-\sum_{k=1}^{p-1}\f{2k+1}k\bi{2k}k\f pk
\\=&\sum_{k=1}^{p-1}\f{2k+1}{k^2}k\bi{2k}k\sum_{j=0}^{p-1}\bi{-k}{j}\bi{-k-1}j-p\sum_{k=1}^{p-1}\f{2k+1}{k^2}\bi{2k}k\pmod{p^2}.
\endalign$$
By [Su16, Lemma 3.4],
$$k\bi{2k}k\sum_{j=0}^{p-1}\bi{-k}j\bi{-k-1}j\eq p\pmod{p^2}\quad\t{for all}\ k=1,\ldots,p-1.$$
So we have
$$\align\sum_{n=1}^{p-1}\f{S_n}n\eq& p\sum_{k=1}^{p-1}\f{2k+1}{k^2}-p\sum_{k=1}^{p-1}\f{2k+1}{k^2}\bi{2k}k
\\\eq&-\f p2\l(\f p3\r)B_{p-2}\l(\f13\r)\pmod{p^2}.
\endalign$$
This concludes the proof of Theorem 1.4(ii). \qed
\medskip

Now we present a $q$-congruence related to (1.21).

\proclaim{Theorem 3.2} Let $a,b\in\N$, and let $n$ be a positive integer. For each $a'\in\{a,a-1\}$, we have
$$\sum_{k=0}^{n-1}(-1)^{a'k}q^{a'k(k+1)/2-k}[2k+1]_q\M{n-1}k^a\M{n+k}k^b\eq0\pmod{[n]_q}.\tag3.7$$
Therefore
$$\sum_{k=0}^{n-1}(\pm1)^k(2k+1)\bi{n-1}k^a\bi{n+k}k^b\eq0\pmod n.\tag3.8$$
\endproclaim
\Proof. (3.8) follows from (3.7) with $q\to1$. Note that (3.7) is trivial for $n=1$.

Below we assume $n>1$ and want to prove (3.7). In view of (3.5), it suffices to show that the left-hand side of (3.7) is divisible by $\Phi_d(q)$
for any divisor $d>1$ of $n$. Write $n=dm$. By Lemma 3.1,
$$\align &\sum_{k=0}^{n-1}(-1)^{a'k}q^{a'k(k+1)/2-k}[2k+1]_q\M{n-1}k^a\M{n+k}k^b
\\=&\sum_{j=0}^{m-1}\sum_{r=0}^{d-1}(-1)^{a'(jd+r)}q^{a'(jd+r)(jd+r+1)/2-(jd+r)}\([2(jd+r)+1]_q
\\&\times\M{(m-1)d+d-1}{jd+r}^a\M{(m+j)d+r}{jd+r}^b\)
\\\eq&\sum_{j=0}^{m-1}(-1)^{a'jd}q^{a'jd(jd+1)/2}
\\&\times\sum_{r=0}^{d-1}(-1)^{a'r}q^{a'r(r+1)/2-r}[2r+1]_q\bi{m-1}j^a\M{d-1}r^a\bi{m+j}j^b\M{r}{r}^b
\\=&\sum_{j=0}^{m-1}(-1)^{a'jd}q^{a'jd(jd+1)/2}\bi{m-1}j^a\bi{m+j}j^b\
\\&\times \sum_{r=0}^{d-1}(-1)^{a'r}q^{a'r(r+1)/2-r}[2r+1]_q\M{d-1}r^a\pmod{\Phi_d(q)}.
\endalign$$
For each $r=0,\ldots,d-1$, we have
$$\align\M{d-1}r=&\prod_{0<s\ls r}\f{1-q^{d-s}}{1-q^s}=\prod_{0<s\ls r}\l(q^{-s}\ \f{q^s-1+(1-q^d)}{1-q^s}\r)
\\\eq&(-1)^rq^{-r(r+1)/2}\pmod{\Phi_d(q)}.\endalign$$
So, by the above, it suffices to show that
$$\sum_{r=0}^{d-1}(-1)^{a'r}q^{a'r(r+1)/2-r}[2r+1]_q\l((-1)^rq^{-r(r+1)/2}\r)^a\eq0\pmod{\Phi_d(q)}.$$
As $a'\in\{a,a-1\}$, this reduces to
$$\sum_{r=0}^{d-1}q^{-r}[2r+1]_q\eq0\eq\sum_{r=0}^{d-1}(-1)^rq^{-r(r+1)/2-r}[2r+1]_q\pmod{\Phi_d(q)}.\tag3.9$$
It is clear that
$$\sum_{r=0}^{d-1}q^{-r}[2r+1]_q=\sum_{r=0}^{d-1}q^{-r}\f{1-q^{2r+1}}{1-q}\eq\sum_{r=0}^{d-1}\f{q^{d-r}-q^{r+1}}{1-q}=0\pmod{\Phi_d(q)}.$$
Also,
$$\align&\sum_{r=0}^{d-1}(-1)^rq^{-(r^2+3r)/2}\f{1-q^{2r+1}}{1-q}
\\=&\f1{1-q}\sum_{r=0}^{d-1}(-1)^r\l(q^{-r(r+3)/2}-q^{-(r-2)(r+1)/2}\r)
\\=&\f1{1-q}\(\sum_{r=0}^{d-1}(-1)^rq^{-r(r+3)/2}-\sum_{r=-2}^{d-3}(-1)^rq^{-r(r+3)/2}\)
\\=&\f1{1-q}\l((-1)^{d-1}q^{-(d-1)(d+2)/2}+(-1)^{d-2}q^{-(d-2)(d+1)/2}\r)
\\=&\f{(-1)^{d-1}}{1-q}\l(q^{1-d(d+1)/2}-q^{1-d(d-1)/2}\r)=(-1)^{d-1}q^{1-d(d+1)/2}[d]_q
\endalign$$
and hence the second congruence in (3.9) holds too.
This concludes the proof. \qed

\heading{4. Proofs of Theorem 1.5 and Corollary 1.1 and some extensions}\endheading

\proclaim{Theorem 4.1} Let $a_1,\ldots,a_m\in\Z$ and $b_1,\ldots,b_m\in\N$. Let $f:\N\to\Z$ be a function with $k\mid f(k)$ for all $k\in\N$.
Let $n$ be a positive integer and set $d=\gcd(a_1,\ldots,a_m,b_1,\ldots,b_m,n)$. Then we have
$$\sum_{k=0}^{n-1}\bar f(k)\prod_{i=1}^m\bi{a_i-1}{b_i+k}\eq0\pmod{d},\tag4.1$$
where
$\bar f(k) = f(k+1)-(-1)^m f(k)$.
If $k^2\mid f(k)$ for all $k\in\N$, then
$$\sum_{k=0}^{n-1}\bar f(k)\prod_{i=1}^m\bi{a_i-1}{b_i+k}\eq(-1)^m\(\sum_{i=1}^ma_i\)\sum_{0<k<n}\f{f(k)}k\prod_{i=1}^m\bi{a_i-1}{b_i+k}\pmod {d^2}.\tag4.2$$
\endproclaim
\Proof. Clearly $f(0)=0$. Observe that
$$\align &\sum_{k=0}^{n-1}\bar f(k)\prod_{i=1}^m\bi{a_i-1}{b_i+k}
\\=&\sum_{k=0}^{n-1}f(k+1)\prod_{i=1}^m\bi{a_i-1}{b_i+k}-(-1)^{m}\sum_{k=0}^{n-1}f(k)\prod_{i=1}^m\bi{a_i-1}{b_i+k}
\\=&\sum_{k=1}^nf(k)\prod_{i=1}^m\bi{a_i-1}{b_i+k-1}-(-1)^{m}\sum_{k=0}^{n-1}f(k)\prod_{i=1}^m\bi{a_i-1}{b_i+k}
\\=&f(n)\prod_{i=1}^m\bi{a_i-1}{b_i+n-1}+\sum_{0<k<n}f(k)d_k,
\endalign$$
where
$$d_k:=\prod_{i=1}^m\l(\bi {a_i}{b_i+k}-\bi{a_i-1}{b_i+k}\r)-(-1)^m\prod_{i=1}^m\bi{a_i-1}{b_i+k}$$
can be written as $\sum_{i=1}^mc_{i,k}\bi{a_i}{b_i+k}$ with $c_{i,k}\in\Z$. Since $k\mid f(k)$ and
$$k\bi {a_i}{b_i+k}= a_i \bi{a_i-1}{b_i+k-1}-b_i\bi{a_i}{b_i+k}\eq0\pmod{d}\tag4.3$$
for all $k=1,2,3,\ldots$, we derive (4.1) from the above.

Now we assume $k^2\mid f(k)$ for all $k\in\N$. For any $0<k<n$, if $1\ls i<j\ls m$ then
$$f(k)\bi{a_i}{b_i+k}\bi{a_j}{b_j+k}=\f{f(k)}{k^2}\(k\bi{a_i}{b_i+k}\)\(k\bi{a_j}{b_j+k}\)\eq0\pmod{d^2},$$
thus we may use (4.3) to deduce that
$$\align f(k)d_k\eq& f(k)\sum_{i=1}^m\bi{a_i}{b_i+k}\prod_{j\not=i}\l(-\bi{a_j-1}{b_j+k}\r)
\\=&\f{f(k)}{k}\sum_{i=1}^m\l(a_i\bi{a_i-1}{b_i+k-1}-b_i\bi{a_i}{b_i+k}\r)(-1)^{m-1}\prod_{j\not=i}\bi{a_j-1}{b_j+k}
\\=&\f{f(k)}{k^2}\sum_{i=1}^m\l(-ka_i\bi{a_i-1}{b_i+k}+(a_i-b_i)k\bi{a_i}{b_i+k}\r)(-1)^{m-1}\prod_{j\not=i}\bi{a_j-1}{b_j+k}
\\\eq&\f{f(k)}k(a_1+\cdots+a_m)(-1)^m\prod_{i=1}^m\bi{a_i-1}{b_i+k}\pmod{d^2}.
\endalign$$
Therefore, (4.2) follows. \qed

\proclaim{Corollary 4.1} Let $a_1,\ldots,a_m\in\Z$ and $b_1,\ldots,b_m\in\N$. Let $n$ be any positive integer and set $d=\gcd(a_1,\ldots,a_m,b_1,\ldots,b_m,n)$. Then we have
$$\align\sum_{k=0}^{n-1}(-1)^{km}\prod_{i=1}^m\bi{a_i-1}{b_i+k}\eq&0\pmod {d},\tag4.4
\\\sum_{k=0}^{n-1}(\pm1)^k(2k+1)\prod_{i=1}^m\bi{a_i-1}{b_i+k}\eq&0\pmod {d},\tag4.5
\\\sum_{k=0}^{n-1}(\pm1)^k(4k^3-1)\prod_{i=1}^m\bi{a_i-1}{b_i+k}\eq&0\pmod {d}.\tag4.6
\endalign$$
Also,
$$\gcd\l(\f{a_1+\cdots+a_m}d-1,2\r)\sum_{k=0}^{n-1}(-1)^{km}(2k+1)\prod_{i=1}^m\bi{a_i-1}{b_i+k}\eq0\pmod {d^2}\tag4.7$$
and
$$6\sum_{k=0}^{n-1}(-1)^{km}(3k^2+3k+1)\prod_{i=1}^m\bi{a_i-1}{b_i+k}\eq0\pmod {d^2}.\tag4.8$$
\endproclaim
\Proof. Clearly, $(-1)^{km+m}=(-1)^{(k+1)m}(k+1)-(-1)^m((-1)^{km}k)$,
$$\align (\pm1)^k(2k+1)=&(\pm1)^{(k+1)-1}(k+1)\pm(\pm1)^{k-1}k,
\\=&(\pm1)^{(k+1)-1}(k+1)^2\mp(\pm1)^{k-1}k^2,
\endalign$$
and
$$\align (\pm1)^k(4k^3-1)=&(\pm1)^{(k+1)-1}(k+1)^2(2(k+1)-3)\pm(\pm1)^{k-1}k^2(2k-3)
\\=&(-1)^{(k+1)-1}\l((k+1)^2k^2-(k+1)\r)\mp(\pm1)^{k-1}\l(k^2(k-1)^2-k\r).
\endalign$$
So (4.4)-(4.6) follow from the first assertion in Theorem 4.1.

Now we prove (4.7). Let $f(k)=(-1)^{km}k^2$ for all $k\in\N$. Then
$$f(k+1)-(-1)^mf(k)=(-1)^{(k+1)m}(2k+1).$$
Applying the second assertion in Theorem 4.1, we get
$$\align&\sum_{k=0}^{n-1}(-1)^{km+m}(2k+1)\prod_{i=1}^m\bi{a_i-1}{b_i+k}
\\\eq&(-1)^m(a_1+\cdots+a_m)\sum_{k=0}^{n-1}(-1)^{km}k\prod_{i=1}^m\bi{a_i-1}{b_i+k}\pmod{d^2}
\endalign$$
and hence
$$\align&\gcd\l(\f{a_1+\cdots+a_m}d-1,2\r)\sum_{k=0}^{n-1}(-1)^{km}(2k+1)\prod_{i=1}^m\bi{a_i-1}{b_i+k}
\\\eq&\f{(a_1+\cdots+a_m)/d}{\gcd((a_1+\cdots+a_m)/d,2)}d\sum_{k=0}^{n-1}(-1)^{km}((2k+1)-1)\prod_{i=1}^m\bi{a_i-1}{b_i+k}\pmod{d^2}.
\endalign$$
Combining this with (4.4) and (4.5), we immediately obtain the desired (4.7).

It remains to show (4.8). Let $g(k)=(-1)^{km}k^3$ for all $k\in\N$. Then
$$g(k+1)-(-1)^mg(k)=(-1)^{(k+1)m}(3k^2+3k+1).$$
Applying the second assertion in Theorem 4.1, we obtain
$$\align&\sum_{k=0}^{n-1}(-1)^{km+m}(3k^2+3k+1)\prod_{i=1}^m\bi{a_i-1}{b_i+k}
\\\eq&(-1)^m(a_1+\cdots+a_m)\sum_{k=0}^{n-1}(-1)^{km}k^2\prod_{i=1}^m\bi{a_i-1}{b_i+k}\pmod{d^2}
\\\eq&0\pmod d
\endalign$$
and hence
$$\align&6\sum_{k=0}^{n-1}(-1)^{km}(3k^2+3k+1)\prod_{i=1}^m\bi{a_i-1}{b_i+k}
\\\eq&\f{a_1+\cdots+a_m}d d\sum_{k=0}^{n-1}(-1)^{km}(2(3k^2+3k+1)-3(2k+1)+1)\prod_{i=1}^m\bi{a_i-1}{b_i+k}
\\\eq&0\pmod{d^2}
\endalign$$
with the use of (4.4) and (4.5). Thus (4.8) holds.

The proof of Corollary 4.1 is now complete. \qed

\Remark\ 4.1. (4.4) was first established by Guo and Zeng [GZ, Theorem 5.5] via $q$-binomial coefficients,
while (4.5) was conjectured by them in [GZ, Conjecture 5.8].

\proclaim{Theorem 4.2} Let $a_1,\ldots,a_m\in\Z$, and let $f:\N\to\Z$ be a function with $k^3\mid f(k)$ for all $k\in\N$.
Then, for any positive integer $n$, we have
$$\aligned&\sum_{k=0}^{n-1}\Delta f(k)\prod_{i=1}^m\bi{a_in-1}{k}\bi{-a_in-1}k
\\\eq& n^2(a_1^2+\cdots+a_m^2)\sum_{0<k<n}\f{f(k)}{k^2}\prod_{i=1}^m\bi{a_in-1}k\bi{-a_in-1}k\pmod {n^3},
\endaligned\tag4.9$$
where $\Delta f(k)=f(k+1)-f(k)$.
\endproclaim
\Proof. Note that $f(0)=0$ and
$$\align &\sum_{k=0}^{n-1}(f(k+1)-f(k)) \prod_{i=1}^m\bi{a_in-1}{k}\bi{-a_in-1}k
\\=&\sum_{k=1}^nf(k)\prod_{i=1}^m\bi{a_in-1}{k-1}\bi{-a_in-1}{k-1}-\sum_{k=0}^{n-1}f(k)\prod_{i=1}^m\bi{a_in-1}{k}\bi{-a_in-1}{k}
\\=&f(n)\prod_{i=1}^m\bi{a_in-1}{n-1}\bi{-a_in-1}{n-1}+\sum_{0<k<n}f(k)d_k(n)-f(0),
\endalign$$
where
$$d_k(n):=\prod_{i=1}^m\bi{a_in-1}{k-1}\bi{-a_in-1}{k-1}-\prod_{i=1}^m\bi{a_in-1}{k}\bi{-a_in-1}{k}.$$
Since
$$\align&\bi{a_in-1}k\bi{-a_in-1}k-\bi{a_in}k\bi{-a_in}k
\\=&\f{a_in-k}k\bi{a_in-1}{k-1}\f{-a_in-k}k\bi{-a_in-1}{k-1}-\f{a_in}k\bi{a_in-1}{k-1}\f{-a_in}k\bi{-a_in-1}{k-1}
\\=&\l(\f{k^2-(a_in)^2}{k^2}+\f{(a_in)^2}{k^2}\r)\bi{a_in-1}{k-1}\bi{-a_in-1}{k-1}=\bi{a_in-1}{k-1}\bi{-a_in-1}{k-1}
\endalign$$
and
$$k^3\bi{a_in}k\bi{-a_in}k\bi{a_jn}k=(a_in)(-a_in)a_jn\bi{a_in-1}{k-1}\bi{-a_in}{k-1}\bi{a_in}{k-1},$$
for $0<k<n$ we have
$$\align k^3d_k(n)=&k^3\prod_{i=1}^m\l(\bi{a_in-1}k\bi{-a_in-1}k-\bi{a_in}k\bi{-a_in}k\r)
\\&-k^3\prod_{i=1}^m\bi{a_in-1}k\bi{-a_in-1}k
\\\eq&-k^3\sum_{i=1}^m\bi{a_in}k\bi{-a_in}k\prod_{j\not=i}\bi{a_jn-1}k\bi{-a_jn-1}k
\\=&n^2\sum_{i=1}^ma_i^2k\l(\bi{a_in}k-\bi{a_in-1}{k}\r)\l(\bi{-a_in}k-\bi{-a_in-1}{k}\r)
\\&\qquad\times\prod_{j\not=i}\bi{a_jn-1}k\bi{-a_jn-1}k
\\\eq&n^2(a_1^2+\cdots+a_m^2)k\prod_{i=1}^m\bi{a_in-1}k\bi{-a_in-1}k\pmod{n^3}.
\endalign$$
Therefore (4.9) follows from the above. \qed

\proclaim{Lemma 4.1} For any $k,n\in\N$, we have
$$\f{k}{\bi{2k-1}k}\bi nk\bi{-n}k\eq0\pmod n.\tag4.10$$
\endproclaim
\Proof. The assertion holds trivially for $k=0$, below we assume $k>0$. In view of (2.6),
$$(-1)^k\bi{n}k\bi{-n}k=\bi{2k-1}k\f{2n}{n+k}\bi{n+k}{2k}=\bi{2k-1}k\f nk\bi{n+k-1}{2k-1}$$
and thus (4.10) follows. \qed

\proclaim{Theorem 4.3} Let $a_1,\ldots,a_m$ be positive integers with $\min\{a_1,\ldots,a_m\}=1$, and let $f$ be a function from $\N$ to the field $\Q$ of rational numbers. Let $n$ be any positive integer.

{\rm (i)} If $\bi{2k-1}kf(k)\in\Z$ for all $k\in\N$, then we have
$$\sum_{k=0}^{n-1}\Delta f(k)\prod_{i=1}^m\bi{a_in-1}{k}\bi{-a_in-1}k\in\Z.\tag4.11$$

{\rm (ii)} If $\bi{2k-1}kf(k)\in k\Z$ for all $k\in\N$, then we have
$$\f1n\sum_{k=0}^{n-1}\Delta f(k)\prod_{i=1}^m\bi{a_in-1}{k}\bi{-a_in-1}k\in\Z.\tag4.12$$
\endproclaim
\Proof. As in the proof of Theorem 4.2, by Abel's partial summation we have
$$\aligned &\sum_{k=0}^{n-1}\Delta f(k)\prod_{i=1}^m\bi{a_in-1}{k}\bi{-a_in-1}k
\\=&f(n)\prod_{i=1}^m\bi{a_in-1}{n-1}\bi{-a_in-1}{n-1}+\sum_{0<k<n}f(k)d_k(n)-f(0),
\endaligned\tag4.13$$
where $$\align d_k(n):=&\prod_{i=1}^m\l(\bi{a_in-1}k\bi{-a_in-1}k-\bi{a_in}k\bi{-a_in}k\r)
\\&-\prod_{i=1}^m\bi{a_in-1}k\bi{-a_in-1}k
\endalign$$ can be written as $\sum_{i=1}^m\bi{a_in}k\bi{-a_in}kc_{i,k}(n)$ with $c_{i,k}(n)\in\Z$.

(i) By Lemma 2.3, $\bi{2k-1}k\mid\bi{a_in}{k}\bi{-a_in}k$ for any $i=1,\ldots,m$ and $k=0,\ldots,n$.
If $f(k)\bi{2k-1}k\in\Z$ for all $k\in\N$, then
$$f(0)\in\Z,\ \ \ f(n)\bi{-n-1}{n-1}=f(n)(-1)^{n-1}\bi{2n-1}n\in\Z,$$
and $f(k)d_k(n)\in\Z$ for all $0<k<n$, thus (4.11) follows (4.13).

(ii) By Lemma 4.1, for any $i=1,\ldots,m$ and $k=0,\ldots,n$ we have
$$\f k{\bi{2k-1}k}\bi{a_in}{k}\bi{-a_in}k\eq0\pmod n.$$
If $\bi{2k-1}kf(k)\in k\Z$ for all $k\in\N$, then $f(0)=0$,
$$(-1)^{n-1}f(n)\bi{-n-1}{n-1}=f(n)\bi{2n-1}n\eq0\pmod n,$$
and $f(k)d_k(n)\eq0\pmod n$ for all $0<k<n$, therefore (4.12) follows from (4.13).

The proof of Theorem 4.3 is now complete. \qed

\proclaim{Theorem 4.4} Let $a$, $b$ and $n$ be positive integers. For any function $f:\N\to\Q$ with $f(k)\bi{2k-1}k\in\Z$ for all $k\in\N$, we have
$$\sum_{k=0}^{n-1}(f(k+1)-(-1)^{a+b}f(k))\bi{n-1}k^a\bi{-n-1}k^b\in\Z.\tag4.14$$
\endproclaim
\Proof. Clearly Theorem 4.3(i) implies (4.14) in the case $a=b$. To handle the general case, we need some new ideas.

By Abel's partial summation,
$$\align&\sum_{k=0}^{n-1}(f(k+1)-(-1)^{a+b}f(k))\bi{n-1}k^a\bi{-n-1}k^b
\\=&\sum_{k=1}^nf(k)\bi{n-1}{k-1}^a\bi{-n-1}{k-1}^b-(-1)^{a+b}\sum_{k=0}^{n-1}f(k)\bi{n-1}k^a\bi{-n-1}k^b
\\=&f(n)\bi{-n-1}{n-1}^b+\sum_{k=0}^{n-1}f(k)\l(\bi nk-\bi{n-1}k\r)^a\l(\bi{-n}k-\bi{-n-1}k\r)^b
\\&-(-1)^{a+b}\sum_{k=0}^{n-1}f(k)\bi{n-1}k^a\bi{-n-1}k^b.
\endalign$$
Note that $\bi{-n-1}{n-1}=(-1)^{n-1}\bi{2n-1}n$. For each $k=0,\ldots,n-1$, we have $\bi{2k-1}k\mid\bi nk\bi{-n}k$ by (2.6), and
$$\align\bi {\pm n}k\bi{\mp n-1}k=&(-1)^k\bi {\pm n}k\bi{\pm n+k}k
\\=&(-1)^k\bi{\pm n+k}{2k}\bi{2k}k=(-1)^k2\bi{\pm n+k}{2k}\bi{2k-1}k,
\endalign$$
therefore
$$\l(\bi nk-\bi{n-1}k\r)^a\l(\bi{-n}k-\bi{-n-1}k\r)^b-(-1)^{a+b}\bi{n-1}k^a\bi{n-1}k^b$$
is divisible by $\bi{2k-1}k$. As $f(k)\bi{2k-1}k\in\Z$ for all $k=0,\ldots,n$, combining the above we obtain (4.14). \qed

\medskip
\noindent{\it Proof of Theorem} 1.5. (i) (1.21)-(1.24) are special cases of (4.5)-(4.8) respectively.
For the function $f(k)=(-1)^{k-1}k^2(2k-3)$, we clearly have $\Delta f(k) = (-1)^k(4k^3-1)$ for all $k\in\N$.
So, (1.25) follows from the last part of Theorem 4.1. As $3k^2+3k+1=(k+1)^3-k^3$, Theorem 4.2 implies that
$$\align &\sum_{k=0}^{n-1}(3k^2+3k+1)\prod_{i=1}^m\bi{a_in-1}k\bi{-a_in-1}k
\\\eq&n^2(a_1^2+\ldots+a_n^2)\sum_{k=0}^{n-1}k\prod_{i=1}^m\bi{a_in-1}k\bi{-a_in-1}k\pmod{n^3}.
\endalign$$
By Corollary 4.1,
$$\sum_{k=0}^{n-1}((2k+1)-1)\prod_{i=1}^m\bi{a_in-1}k\bi{-a_in-1}k\eq0\pmod n.$$
Therefore
$$\align&\gcd(a_1+\cdots+a_m-1,2)\sum_{k=0}^{n-1}(3k^2+3k+1)\prod_{i=1}^m\bi{a_in-1}k\bi{-a_in-1}k
\\\eq&n^2\f{a_1^2+\cdots a_m^2}{\gcd(a_1+\cdots+a_m,2)}\sum_{k=0}^{n-1}((2k+1)-1)\prod_{i=1}^m\bi{a_in-1}k\bi{-a_in-1}k
\\\eq&0\pmod{n^3}.
\endalign$$
This proves (1.26).

(ii) Now let $a,b,n$ be positive integers. Note that
 $$\f2{k+1}\bi{2k-1}k=\f{\bi{2k}k}{k+1}=C_k
 \ \ \t{and}\ \ \f{\bi{2k-1}k}{2k-1}=\cases C_{k-1}&\t{if}\ k>0,\\-1&\t{if}\ k=0.\endcases$$
For $k\in\N$, define
 $$f_1(k)=\f k{2k-1},\ f_2(k)=\f{(-1)^kk}{2k-1},\ f_3(x)=\f {2k}{k+1},\ f_4(x)=\f{(-1)^k2k}{k+1}.$$
Then $f_i(k)\bi{2k-1}k\in k\Z$ for all $i=1,\ldots,4$. Clearly,
 $$\align\Delta f_1(k)=&\f {k+1}{2k+1}-\f k{2k-1}=-\f1{4k^2-1},
 \\ \Delta f_2(k)=&\f{(-1)^{k+1}(k+1)}{2k+1}-\f{(-1)^kk}{2k-1}=(-1)^{k-1}\l(1+\f{2k}{4k^2-1}\r),
 \\\Delta f_3(k)=&\f{2(k+1)}{k+2}-\f{2k}{k+1}=\f1{\bi{k+2}2},
 \\ \Delta f_4(k)=&\f{(-1)^{k+1}2(k+1)}{k+2}-\f{(-1)^k2k}{k+1}=(-1)^{k-1}\l(4-\f{2k+3}{\bi{k+2}2}\r).
 \endalign$$
 Applying Theorem 4.3(ii) with $f=f_1,\ldots,f_4$, we immediately get (1.27)-(1.29).

  Write $m=a+b$. For $k\in\N$, define
 $$\gather f_5(k)=\f{(-1)^{km}}{2k-1},\ f_6(k)=\f{(-1)^{k(m-1)}}{2k-1},\ f_7(k)=\f {(-1)^{km}2}{k+1},
 \\\ f_8(k)=\f{(-1)^{k(m-1)}2}{k+1},\ f_9(k)=\f{(-1)^{km}}{\bi{2k-1}k},\ f_{10}(k)=\f{(-1)^{k(m-1)}}{\bi{2k-1}k}.\endgather$$
 Then $f_i(k)\bi{2k-1}k\in\Z$ for all $i=5,\ldots,10$.
 Let $\bar f_i(k)=f_i(k+1)-(-1)^mf_i(k)$ for $i=5,\ldots,10$. Observe that
 $$\align \bar f_5(k)=&\f{(-1)^{(k+1)m}}{2k+1}-(-1)^m\f{(-1)^{km}}{2k-1}=(-1)^{(k-1)m}\f{-2}{4k^2-1},
 \\\bar f_6(k)=&\f{(-1)^{(k+1)(m-1)}}{2k+1}-(-1)^m\f{(-1)^{k(m-1)}}{2k-1}=(-1)^{(k-1)(m-1)}\f{4k}{4k^2-1},
 \\\bar f_7(k)=&\f{(-1)^{(k+1)m}2}{k+2}-(-1)^m\f{(-1)^{km}2}{k+1}=(-1)^{(k-1)m}\f{-1}{\bi{k+2}2},
 \\\bar f_8(k)=&\f{(-1)^{(k+1)(m-1)}2}{k+2}-(-1)^m\f{(-1)^{k(m-1)}2}{k+1}=(-1)^{(k-1)(m-1)}\f{2k+3}{\bi{k+2}2},
 \\\bar f_9(k)=&\f{(-1)^{(k+1)m}}{\bi{2k+1}{k+1}}-(-1)^m\f{(-1)^{km}}{\bi{2k-1}k}=(-1)^{(k-1)m}\f{-(3k+1)}{(2k+1)\bi{2k}k},
 \endalign$$
 and
 $$\bar f_{10}(k)=\f{(-1)^{(k+1)(m-1)}}{\bi {2k+1}{k+1}}-(-1)^m\f{(-1)^{k(m-1)}}{\bi{2k-1}k}=\f{(-1)^{(k-1)(m-1)}(5k+3)}{(2k+1)\bi{2k}k}.$$
Theorem 4.4 with $f=f_5,\ldots,f_{10}$ clearly yields (1.30)-(1.35).

The proof of Theorem 1.5 is now complete. \qed

\proclaim{Lemma 4.2} Let $a_0,a_1,\ldots$ be a sequence of complex numbers, and define
$$\tilde a_n=\sum_{k=0}^n\bi nk^2\bi{n+k}k^2a_k\quad\t{for}\ n\in\N.\tag4.15$$
Then, for any positive integer $n$, we have
$$\f1{n^2}\sum_{k=0}^{n-1}(2k+1)\tilde a_k=\sum_{k=0}^{n-1}\f{a_k}{2k+1}\bi{n-1}k^2\bi{n+k}k^2.\tag4.16$$
\endproclaim
\Proof. By [Su12b, Lemma 2.1],
$$\sum_{m=0}^{n-1}(2m+1)\bi{m+k}{2k}^2=\f{(n-k)^2}{2k+1}\bi{n+k}{2k}^2\quad\t{for all}\ k\in\N.$$
Thus
$$\align\sum_{m=0}^{n-1}(2m+1)\tilde a_m=&\sum_{m=0}^{n-1}(2m+1)\sum_{k=0}^m\bi{m+k}{2k}^2\bi{2k}k^2a_k
\\=&\sum_{k=0}^{n-1}\bi{2k}k^2a_k\sum_{m=0}^{n-1}(2m+1)\bi{m+k}{2k}^2
\\=&\sum_{k=0}^{n-1}\bi{2k}k^2\f{a_k}{2k+1}(n-k)^2\bi{n+k}{2k}^2
\\=&\sum_{k=0}^{n-1}\f{a_k}{2k+1}(n-k)^2\bi{n}k^2\bi{n+k}k^2
\\=&n^2\sum_{k=0}^{n-1}\f{a_k}{2k+1}\bi{n-1}k^2\bi{n+k}k^2.
\endalign$$
This proves (4.16). \qed

\medskip\noindent
{\it Proof of Corollary} 1.1. By Lemma 4.2 and (1.27), we have
$$\f1{n^3}\sum_{k=0}^{n-1}(2k+1)t_k=\f1n\sum_{k=0}^{n-1}\f{\bi{n-1}k^2\bi{n+k}k^2}{4k^2-1}\in\Z.$$
In light of Lemma 4.2,
$$\f1{n^2}\sum_{k=0}^{n-1}(2k+1)T_k=\sum_{k=0}^{n-1}\bi{n-1}k^2\bi{n+k}k^2.$$
By [GZ, (1.9)] or (4.4),
$$\sum_{k=0}^{n-1}\bi{n-1}k^2\bi{n+k}k^2\eq0\pmod n.$$
So we have $\sum_{k=0}^{n-1}(2k+1)T_k\eq0\pmod{n^3}$.
By Lemma 4.2 and (1.36) and (1.21),
$$\f1{n^4}\sum_{k=0}^{n-1}(2k+1)T_k^+=\f1{n^2}\sum_{k=0}^{n-1}(2k+1)\bi{n-1}k^2\bi{n+k}k^2\in\Z$$
and
$$\f1{n^3}\sum_{k=0}^{n-1}(2k+1)T_k^{-}=\f1n\sum_{k=0}^{n-1}(-1)^k(2k+1)\bi{n-1}k^2\bi{n+k}k^2\in\Z.$$
Therefore both (1.37) and (1.38) hold. This concludes the proof. \qed

\heading{5. Some related conjectures}\endheading

\proclaim{Conjecture 5.1} Let $p\eq3\pmod 4$ be a prime. Then
$$\sum_{k=0}^{p-1}\f{\bi{2k}k^2}{(2k-1)8^k}\eq-\l(\f 2p\r)\f{p+1}{2^{p-1}+1}\bi{(p+1)/2}{(p+1)/4}\pmod{p^2}\tag5.1$$
and
$$3\sum_{k=0}^{p-1}\f{\bi{2k}k\bi{2k}{k+1}}{(2k-1)8^k}\eq p+\l(\f 2p\r)\f{2p}{\bi{(p+1)/2}{(p+1)/4}}\pmod{p^2}.\tag5.2$$
\endproclaim

\proclaim{Conjecture 5.2} {\rm (i)} The sequence $(R_{n+1}/R_n)_{n\gs3}$ is strictly increasing to the limit $3+2\sqrt2$, and the sequence
$(\root{n+1}\of{R_{n+1}}/\root n\of{R_n})_{n\gs5}$ is strictly decreasing.

{\rm (ii)} The sequence $(S_{n+1}/S_n)_{n\gs3}$ is strictly increasing to the limit $9$, and the sequence
$(\root{n+1}\of{S_{n+1}}/\root n\of{S_n})_{n\gs1}$ is strictly decreasing.
\endproclaim
\Remark\ 5.1. The author [Su13b] made many similar conjectures for some well-known integer sequences.

\proclaim{Conjecture 5.3} For any positive integer $n$, both $R_n(x)$ and $S_n(x)$ are irreducible over the field of rational numbers.
\endproclaim

\proclaim{Conjecture 5.4} For any $n\in\Z^+$, the number $\f3n\sum_{k=0}^{n-1}R_k^2$ is always an odd integer; moreover,
$$\f 3n\sum_{k=0}^{n-1}R_k(x)^2\in\Z[x]\quad\t{and}\quad\f1n\sum_{k=0}^{n-1}(2k+1)R_k^2\in\Z.\tag5.3$$
Also, for any odd prime $p$ we have
$$\sum_{k=0}^{p-1}R_k^2\eq\f p3\l(11-4\l(\f{-1}p\r)\r)\pmod{p^2}\tag5.4$$
and
$$\sum_{k=0}^{p-1}(2k+1)R_k^2\eq 4p\l(\f{-1}p\r)-p^2\pmod{p^3}.\tag5.5$$
\endproclaim
\Remark\ 5.2. For any positive integer $n$, we can easily deduce that
$$\f3n\sum_{k=0}^{n-1}(2k+1)R_k(x)=\sum_{k=0}^{n-1}(n-k)\bi{n+k}{2k}\bi{2k}k\l(\f2{2k-1}-\f1{k+1}\r)x^k\in\Z[x].\tag5.6$$

\proclaim{Conjecture 5.5} We have
$$\f4{n^2}\sum_{k=0}^{n-1}kS_k\in\Z\quad\t{for all}\ n=1,2,3,\ldots.\tag5.7$$
Also, for any prime $p$ we have
$$\sum_{k=0}^{p-1}kS_k\eq\f{p^2}8\l(5-9\l(\f p3\r)\r)\pmod{p^3}.\tag5.8$$
\endproclaim

\proclaim{Conjecture 5.6} For $n\in\N$ define
$$\align s_n:=&\sum_{k=0}^n\bi nk^2\bi{2k}k\f1{2k-1},
\\S_n^+:=&\sum_{k=0}^n\bi nk^2\bi{2k}k(2k+1)^2,
\\S_n^-:=&\sum_{k=0}^n\bi nk^2\bi{2k}k(2k+1)^2(-1)^k.
\endalign$$
Then, for any positive integer $n$, we have
$$\f1{n^2}\sum_{k=0}^{n-1}s_k\in\Z,\ \ \f1{n^2}\sum_{k=0}^{n-1}S_k^+\in\Z\ \ \t{and}\ \ \f1{n^2}\sum_{k=0}^{n-1}S_k^-\in\Z.\tag5.9$$
\endproclaim
\Remark\ 5.3. For any positive integer $n$, we can easily deduce $\sum_{k=0}^{n-1}S_k^{\pm}\eq0\pmod n$
with the help of (3.4). We also conjecture that $\sum_{k=0}^{p-1}s_k\eq-(9(\f p3)+1)p^2/2\pmod{p^3}$ for any prime $p$.

\proclaim{Conjecture 5.7} For $n\in\N$ define
$$s_n(q):=\sum_{k=0}^n\M nk^2\M{2k}k\f{q^k}{[2k-1]_q}.$$
Then, for any positive integer $n$, we have
$$\f{1+q}2\sum_{k=0}^{n-1}q^ks_k(q)\eq0\pmod{[n]_q^2}.\tag5.10$$
\endproclaim
\Remark\ 5.4. (5.10) is a $q$-analogue of the conjectural congruence $\sum_{k=0}^{n-1}s_k\eq0\pmod {n^2}$. We could prove (5.10) modulo $[n]_q$.

\proclaim{Conjecture 5.8} Let $m$ be any positive integer.

{\rm (i)} Define
$$S^{(m)}_n(x):=\sum_{k=0}^n\bi nk^m\f{(km+1)!}{(k!)^m}x^k\quad\t{for}\ \ n=0,1,2,\ldots.$$
Then, for any positive integer $n$, we have
$$\f1n\sum_{k=0}^{n-1}S^{(m)}_k(x)\in\Z[x],\tag5.11$$
i.e.,
$$\f{(km+1)!}{(k!)^m}\sum_{h=k}^{n-1}\bi hk^m\eq0\pmod n\quad\t{for all}\ \ k=0,\ldots,n-1.\tag5.12$$

{\rm (ii)} Define
$$S^{(m)}_n(x;q)=\sum_{k=0}^n\M nk^m\f{\prod_{j=1}^{km+1}[j]_q}{(\prod_{0<j\ls k}[j]_q)^m}x^k\quad\t{for}\ \ n=0,1,2,\ldots.\tag5.13$$
Then, for any integer $n>0$, all the coefficients of the polynomial $\sum_{k=0}^{n-1}q^kS^{(m)}_k(x;q)$ in $x$ are divisible by $[n]_q$ in the ring $\Z[q]$, i.e.,
$$\f{\prod_{j=1}^{km+1}[j]_q}{(\prod_{0<j\ls k}[j]_q)^m}\sum_{h=k}^{n-1}q^h\M hk^m\eq0\pmod {[n]_q}\quad\t{for all}\ k=0,\ldots,n-1.\tag5.14$$
\endproclaim
\Remark\ 5.5. (a) Note that $S_n^{(2)}(x)=S_n(x)$, and (5.11) and (5.12) are extensions of (1.20) and (3.4) respectively. Part (ii) of Conjecture 5.8
presents a $q$-analogue of the first part, and our Theorem 3.1 confirms it for $m=2$. Conjecture 5.8 for $m=1$ is easy, and we are also able to prove
Conjecture 5.8 in the case $m=3$.

(b) The congruence in (5.12) for $k=1$ states that
$$(m+1)!\sum_{h=1}^{n-1}h^m\eq0\pmod n.$$
This is easy since
$$\f1n\sum_{h=0}^{n-1}h^m=\f1{m+1}\sum_{k=0}^m\bi{m+1}kB_kn^{m-k}$$
(cf. [IR, p.\,230]) and $(k+1)!B_k\in\Z$ by the von Staudt-Clausen theorem (cf. [IR, p.\,233]).

\Ack. The author would like to thank his graduate student Xiang-Zi Meng and the referee for helpful comments.
The initial version of this paper was posted to arXiv as a preprint with the ID {\tt arXiv:1408.5381}.


\widestnumber\key{Su12b}

 \Refs

\ref\key AAR\by G. E. Andrews, R. Askey and R. Roy\book Special Functions\publ Cambridge Univ. Press, Cambridge, 1999\endref

\ref\key BEW\by B. C. Berndt, R. J. Evans and K. S. Williams
\book Gauss and Jacobi Sums\publ John Wiley \& Sons, 1998\endref

\ref\key CDE\by S. Chowla, B. Dwork and R. J. Evans\paper On the mod $p^2$ determination of $\bi{(p-1)/2}{(p-1)/4}$
\jour J. Number Theory\vol24\yr 1986\pages 188--196\endref

\ref\key G\by H.W. Gould\book {\rm Combinatorial Identities} \publ
Morgantown Printing and Binding Co., 1972\endref

\ref\key GZ\by V.J.W. Guo and J. Zeng\paper New congruences for sums involving Ap\'ery numbers or central Delannoy numbers
\jour Int. J. Number Theory \vol 8\yr 2012\pages 2003--2016\endref

\ref\key IR\by K. Ireland and M. Rosen \book {\rm A Classical
Introduction to Modern Number Theory (Graduate Texts in
Math.; 84), 2nd ed.} \publ Springer, New York, 1990\endref

\ref\key MT\by S. Mattarei and R. Tauraso\paper Congruences for central binomial sums and finite polylogarithms
\jour J. Number Theory\vol 133\yr 2013\pages 131--157\endref


\ref\key O\by G. Olive\paper Generalized powers\jour Amer. Math. Monthly \vol 72 \yr 1965\pages 619--627\endref

\ref\key PWZ\by M. Petkov\v sek, H. S. Wilf and D. Zeilberger\book
$A=B$ \publ A K Peters, Wellesley, 1996\endref


\ref\key S11\by Z.-H. Sun\paper Congruences concerning Legendre
polynomials \jour Proc. Amer. Math. Soc. \vol 139\yr 2011\pages 1915--1929\endref

\ref\key Su11\by Z.-W. Sun\paper On congruences related to central
binomial coefficients \jour J. Number Theory\vol 131\yr 2011\pages
2219--2238\endref

\ref\key Su12a\by Z.-W. Sun\paper On sums involving products of three
binomial coefficients \jour Acta Arith. \vol 156\yr 2012\pages 123--141\endref

\ref\key Su12b\by Z.-W. Sun\paper On sums of Ap\'ery polynomials and related congruences
 \jour J. Number Theory \vol 132\yr 2012\pages 2673--2699\endref

\ref\key Su13a\by Z.-W. Sun\paper Supercongruences involving products of two binomial coefficients
\jour Finite Fields Appl.\vol 22\yr 2013\pages 24--44\endref

\ref\key Su13b\by Z.-W. Sun\paper Conjectures involving arithmetical sequences\jour
in: Number Theory: Arithmetic in Shangri-La (Shanghai, 2011), S. Kanemitsu et al. (eds.), World Sci.,
Hackensack, NJ, 2013, pp. 244-258\endref

\ref\key Su16\by Z.-W. Sun\paper Congruences involving $g_n(x)=\sum_{k=0}^n\bi nk^2\bi{2k}kx^k$
\jour Ramanujan J.\vol 40\yr 2016\pages 511--533\endref

\ref\key ST\by Z.-W. Sun and R. Tauraso\paper New congruences for central binomial coefficients
\jour Adv. Appl. Math. \vol 45\yr 2010\pages 125--148\endref

\ref\key W\by J. Wolstenholme\paper On certain properties of prime numbers\jour Quart. J. Appl. Math.
\vol 5\yr 1862\pages 35--39\endref

\endRefs
\enddocument